\documentclass[12pt]{amsart}
\usepackage[cp1251]{inputenc}
\usepackage[english]{babel}
\usepackage{amsmath}
\usepackage{amssymb}
\usepackage{amsfonts}

\begin{document}

\title[Operator symbols. II.]
      {Operator symbols. II.}
\author{Vladimir~B. Vasilyev}
\address{Chair of Applied Mathematics and Computer Modeling,
 Belgorod State National Research University,
Pobedy street 85, Belgorod 308015, Russia}

\email{vbv57@inbox.ru}

\begin{abstract}
 We consider special classes of linear bounded operators in Banach spaces and suggest certain operator variant of symbolic calculus. It permits to formulate an index theorem and to describe Fredholm properties of elliptic pseudo-differential operator on manifolds with non-smooth boundaries.
\end{abstract}

\maketitle

\section{Introduction}

In this paper we consider some abstract operators acting in some functional spaces. These considerations were inspired by studies of I.B. Simonenko \cite{S} related to special operators of a local type (we say here local operators). Such operators and corresponding equations plays important role in the theory of pseudo-differential operators and equations \cite{E,RS,V1}. There are a few approaches to the theory of pseudo-differential operators and equations on non-smooth manifolds and manifolds with non-smooth boundarues, but it seems the suggested abstract variant is very close to this theory. Some first steps were done in the author's preprint \cite{V6}, and here we develop this abstract variant and give some applications. We think this approach can be useful for similar problems related to concrete operators.

\section{Operator symbols}

\subsection{Local operators}
In this section we will give some preliminary ideas and definitions from \cite{S,V6}.  {\bf Here we consider such functional spaces which include smooth functions and corresponding multipliers and only local operators.       Additionally all considered operators are defined up to compact operators.}

Let $M$ be a compact $m$-dimensional manifold with a boundary. Below we will consider the case of piecewise smooth boundary, and all singularities will be described. Here we will try to develop certain general statements.

Let $B_1,B_2$ be  Banach spaces consisting of functions defined on compact $m$-dimensional manifold $M$. We assume
that smooth functions with compact support are dense in such spaces. Let $A: B_1\rightarrow B_2$ be a
linear bounded operator\footnote{We remind that an index of the operator $A$ is called the following number $dim~ Ker A-dim~ Coker A$ \cite{E,RS,V1}}. We will denote by letter $f$ the function $f$ and the operator of multiplication by $f$, so that the notation $A\cdot f$ denotes the following operator
\[
(A\cdot f)u=A(fu),~~~u\in B_1.r
\]

{\bf Definition 1.} {\it
An operator $A$ is called a local operator if the operator
\[
f\cdot A\cdot g
\]
is a compact operator for arbitrary smooth functions with non-intersecting supports.
}

\subsection{Operators on a compact manifold}

On the manifold $M$ we fix a finite open covering and a partitions of unity corresponding to this covering $\{U_j, f_j\}_{j=1}^{n}$ and
choose smooth functions
  $\{ g_j\}_{j=1}^{n}$ so that $supp~g_j\subset V_j,$ $\overline{U_j}\subset V_j$, and $g_j(x)\equiv 1$ for $x\in supp~ f_j, supp~f_j\cap(1-g_j)=\emptyset$.

{\bf Proposition 1.} {\it
 The  operator $A$ on the manifold $M$ can be represented in the form
\[
A=\sum\limits_{j=1}^{n}f_j\cdot A\cdot g_j +T,
\]
where $T: B_1\rightarrow B_2$ is a compact operator.
}

{\bf Proof.}
It is very simple. Since
\[
\sum\limits_{j=1}^nf_j\equiv 1,
\]
then we can write
\[
\left(\sum\limits_{j=1}^nf_j\right)\cdot A=\sum\limits_{j=1}^nf_j\cdot A=\sum\limits_{j=1}^nf_j\cdot A\cdot g_j+\sum\limits_{j=1}^nf_j\cdot A\cdot(1-g_j),
\]
so we have the conclusion needed.
$\triangle$

{\bf Remark 1.} {\it
It is obvious $A$ is defined uniquely up to a compact operators which do not influence on an index.
}

By definition for an arbitrary operator $A: B_1\rightarrow B_2$
\[
|||A|||\equiv\inf||A+T||,
\]
where {\it infimum} is taken over all compact operators $T: B_1\rightarrow B_2$.

Let $B'_1,B'_2$ be  Banach spaces consisting of functions defined on ${\bf R}^m$, $\widetilde A: B'_1\rightarrow B'_2$ be a linear bounded operator.

Since $M$ is a compact manifold, then for every point $x\in M$ there exists a neighborhood $U\ni x$ and diffeomorphism $\omega: U\rightarrow D\subset{\bf R}^m, \omega(x)\equiv y$. We denote by $S_{\omega}$ the following operator\footnote{Really, this operator is defined locally; in general it may be unbounded $B_k\rightarrow B'_k$ (see \cite{V5})} acting from $B_k$ to $B_k', k=1,2$. For every function $u\in B_k$ vanishing out of $U$
\[
(S_{\omega}u)(y)=u(\omega^{-1}(y)),~~~y\in D,~~~(S_{\omega}u)(y)=0,~~~y\notin D.
\]

Of course, for every function $v\in B_k'$ vanishing out of $D$ we can define
\[
(S^{-1}_{\omega}v)(x)=v(\omega(x)),~~~x\in U,~~~(S^{-1}_{\omega}v)(x)-0,~~~x\notin U.
\]

{\bf Definition 2.} {\it
A local representative of the operator $A: B_1\rightarrow B_2$ at the point $x\in M$  is called the operator $\widetilde A: B_1'\rightarrow B_2'$ such that $\forall\varepsilon >0$ there exists the neighborhood $U_j$ of the point $x\in  U_j\subset M$ and diffeomorphism $\omega_j': U_j\rightarrow D_j\subset{\bf R}^m$ with the property
\[
|||g_jAf_j-S^{-1}_{\omega_j}\hat g_j\widetilde A\hat f_jS_{\omega_j}|||<\varepsilon.
\]
}

\section{Generating operator}

      Let $M$ be a compact $m$-dimensional manifold  with a boundary $\partial M$, and $A(x)$ be a certain operator-function defined on $M$. Let $M_k, k=0,1,...,m-1,$ be smooth $k$-dimensional sub-manifolds on $\partial M$ so that by definition $M_{m-1}\equiv\partial M, M_0$ consists of isolated points on $\partial M$. Further, we introduce a set of operator classes ${\rm T}_k, k=0,1,...,m,$ so that for $x\in M_k$, $A(x): H^{(1)}_k\rightarrow  H^{(2)}_k$ is a linear bounded operator, where
$H^{(j)}_k, k=0,1,...,m, j=1,2,$ are some Banach spaces.

We say that sub-manifold $M_k$ is a singular $k$-sub-manifold if $\forall x\in M_k$ we have $A(x)\in{\rm T}_k$. Additionally, we will assume that if $x\in M_r\cap M_{k-1}\neq\emptyset$ then $A(x)\in{\rm T}_{k-1}$.

{\bf Theorem 1.} {\it
If the family $A(x)$ consists of local Fredholm operators and this family is continuous on each component $\overline{M_k\setminus\cup_{i=0}^{k-1}M_i}, k=0,1,...,m,$ then it generates a unique Fredholm operator
$A$ acting in the spaces $\sum\limits_{k=0}^m\oplus H^{(1)}_k\rightarrow\sum\limits_{k=0}^m\oplus H^{(2)}_k$.
}

{\bf Proof.}
First, we construct such an operator in the following way. Let $\varepsilon>0$ is enough small. We take a covering for $M$ by balls as follows. We take a covering for $M_0$, it consists of finite number of open sets and denote this covering by $\mathcal U_0$. Further, we compose $M\setminus\mathcal U_0$. For every point $x\in M_1\cap(M\setminus\mathcal U_0)$ we take ball with the center $x$ of radius $\varepsilon$. The union of such balls is covering for the set $x\in M_1\cap(M\setminus\mathcal U_0)$. According to compactness of the set we extract a finite sub-covering which will be denoted by $\mathcal U_1$. Then we compose the set $M_2\cap( M\setminus(\mathcal U_0\cup\mathcal U_1))$, repeate the procedure mentioned above and obtain the sub-covering $\mathcal U_2$. Continuing the process we obtain the finite covering for $M$ of the following type
\[
M\subset\bigcup\limits_{k=0}^m\mathcal U_k\equiv\mathcal U
\]
Without loss of generality we can mean that elements of the covering are balls with centers at points $x^{(k)}_{j}\in M_k, j=0,1,\cdots,n_k, k=0,1,\cdots,m$.

Since the set $M_0$ consists of isolated points only we have a finite number of operators acting $H^{(1)}_0\rightarrow  H^{(2)}_0$. We construct a partition of unity $f^{(k)}_{j}$ for every sub-covering $\mathcal U_k$ and associated set of functions $g^{(k)}_{j}, j=0,1,\cdots,n_k, k=1,2,\cdots,m$. Let us consider the $k$th component.

Using a piece of the operator-function $A(x)$ related to $M_k$ we construct the following sequence of operators acting $H^{(1)}_k\rightarrow  H^{(2)}_k$. Let us denote
\[
A_{n_k}=\sum\limits_{j-1}^{n_k}f^{(k)}_{j}\cdot A(x^{(k)}_j)\cdot g^{(k)}_{j}
\]
and consider another sub-covering $\mathcal V_k$ for the set $M\setminus(\bigcup\limits_{l=0}^{k-1}\mathcal U_l)$.
Let is suppose that this covering consists of balls with centers in $y^{(k)}_{i}\in M_k, i=1,2,\cdots,r_k$ of enough small radius. We can construct the operator
\[
A_{r_k}=\sum\limits_{i-1}^{r_k}f^{(k)}_{i}\cdot A(y^{(k)}_i)\cdot g^{(k)}_{i}
\]

We would like to prove the following sentence
\begin{equation}\label{4}
|||A_{n_k}-A_{r_k}|||\to 0, ~~~\text{if}~~~ n_k, r_k\to\infty.
\end{equation}
under appropriate choice of coverings $\mathcal U_k, \mathcal V_k$.

As soon as the formula \eqref{4} will be proved we conclude that the sequence $\{A_{n_k}\}$ is a Cauchy sequence with respect to the norm $|||\cdot|||$. Therefore there exists the operator limit $A^{(k)}=\lim\limits_{n_k\to\infty}A_{n_k}$.

The left part of the proof repeats, in general, arguments from \cite{V6}, but for reader's convenience we give these reasonings here in view of their values.

We will construct the $k$th component for the operator $A$ in the following way. Let $\{\varepsilon_n\}_{n=1}^{\infty}$  be a sequence such that $\varepsilon_n>0, \forall n\in{\bf N}, \lim\limits_{n\to\infty}\varepsilon_n=0$. Given $\varepsilon_n$ we choose coverings $\{U^{(k)}_{j}\}_{j=1}^{n_k}\equiv\mathcal U_k$  as above with partition of unity $\{f^{(k)}_{j}\}$ and corresponding functions $\{g^{(k)}_{j}\}$ such that
\[
|||f^{(k)}_{j}\cdot(A(x)-A(x^{(k)}_{j}))\cdot g^{(k)}_{j}|||<\varepsilon_{n_k},~~~\forall x\in U^{(k)}_{j},
\]
and $\{V^{(k)}_{i}\}_{i=1}^{r_k}\equiv\mathcal V_k$ with partition of unity $\{F^{(k)}_{i}\}$ and corresponding functions $\{G^{(k)}_{i}\}$ such that
\[
|||F^{(k)}_{i}\cdot(A(x)-A(y^{(k)}_{i}))\cdot G^{(k)}_{i}|||<\varepsilon_{r_k},~~~\forall x\in V^{(k)}_{i};
\]
we remaind that $U^{(k)}_{j}, V^{(k)}_{i}$   are balls with centers at $x^{(k)}_{j}, y^{(k)}_{i}\in\overline{M_k}$ of radius $\varepsilon$ and $2\varepsilon$. This requirement is possible according to continuity of the operator family $A(x)$ with respect to the norm $|||\cdot|||$ on the sub-manifold $\overline{M_k}$.

We can write
\[
A_{n_k}=\sum\limits_{j=1}^{n_k}f^{(k)}_{j}\cdot A(x^{(k)}_{j})\cdot g^{(k)}_{j}=\sum\limits_{i=1}^{r_k}F^{(k)}_{i}\cdot\sum\limits_{j=1}^{n_k}f^{(k)}_{j}\cdot A(x^{(k)}_{j})\cdot g^{(k)}_{j}=
\]
\[
\sum\limits_{i=1}^{r_k}\sum\limits_{j=1}^{n_k}F^{(k)}_{i}\cdot f^{(k)}_{j}\cdot A(x^{(k)}_{j})\cdot g^{(k)}_{j}=\sum\limits_{i=1}^{r_k}\sum\limits_{j=1}^{n_k}F^{(k)}_{i}\cdot f^{(k)}_{j}\cdot A(x^{(k)}_{j})\cdot g^{(k)}_{j}\cdot G^{(k)}_i+T_1,
\]
and the same we can write for $A_{r_k}$
\[
A_{r_k}=\sum\limits_{i=1}^{r_k}F^{(k)}_{i}\cdot A(y^{(k)}_{i})\cdot G^{(k)}_{i}=\sum\limits_{j=1}^{n_k}f^{(k)}_{j}\cdot \sum\limits_{i=1}^{r_k}F^{(k)}_{i}\cdot A(y^{(k)}_{i})\cdot G^{(k)}_{i}=
\]
\[
\sum\limits_{j=1}^{n_k}\sum\limits_{i=1}^{r_k}f^{(k)}_{j}\cdot F^{(k)}_{i}\cdot A(y^{(k)}_{i})\cdot G^{(k)}_{i}=\sum\limits_{j=1}^{n_k}\sum\limits_{i=1}^{r_k}f^{(k)}_{j}\cdot F^{(k)}_{i}\cdot A(y^{(k)}_{i})\cdot G^{(k)}_{i}\cdot g^{(k)}_{j}+T_2.
\]

Let us consider the difference
\begin{equation}\label{6}
|||A_{n_k}-A_{r_k}|||=|||\sum\limits_{j=1}^{n_k}\sum\limits_{i=1}^{r_k}f^{(k)}_{j}\cdot F^{(k)}_{i}\cdot (A(x^{(k)}_{j})-A(y^{(k)}_{i}))\cdot G^{(k)}_{i}\cdot g^{(k)}_{j}|||.
\end{equation}

Obviously, summands with non-vanishing supplements to the formula \eqref{6} are those for which $U^{(k)}_{j}\cap V^{(k)}_{i}\neq\emptyset$. A number of such neighborhoods are finite always for arbitrary finite coverings, hence we obtain
\[
|||A_{n_k}-A_{r_k}|||\leq\sum\limits_{j=1}^{n_k}\sum\limits_{i=1}^{r_k}|||f^{(k)}_{j}\cdot F^{(k)}_{i}\cdot (A(x^{(k)}_{j})-A(y^{(k)}_{i}))\cdot G^{(k)}_{i}\cdot g^{(k)}_{j}|||\leq
\]
\[
\sum\limits_{x\in U^{(k)}_{j}\cap V^{(k)}_{i}\neq\emptyset}|||f^{(k)}_{j}\cdot F^{(k)}_{i}\cdot (A(x^{(k)}_{j})-A(x))\cdot G^{(k)}_{i}\cdot g^{(k)}_{j}|||+
\]
\[
\sum\limits_{x\in U^{(k)}_{j}\cap V^{(k)}_{i}\neq\emptyset}|||f^{(k)}_{j}\cdot F^{(k)}_{i}\cdot (A(x)-A(y^{(k)}_{i}))\cdot G^{(k)}_{i}\cdot g^{(k)}_{j}|||\leq 2K\max[\varepsilon_{n_k},\varepsilon_{r_k}],
\]
where $K$ is a universal constant.

Thus, we have proved that the sequence $\{A_{n_k}\}$ is a Cauchy sequence, hence there exists $\lim\limits_{{n_k}\to\infty}A_{n_k}=A^{(k)}$.

Using the same process we can construct all operators $A^{(k)}$ for every $k=0,1,\cdots,m$. Let us note all operators $A^{(k)}: H^{(1)}_k\rightarrow H^{(2)}_k$ act in different spaces. Finally, it is easy to compose the resulting operator $A$ acting in direct sums of such spaces. Indeed, if
\[
u=\oplus\sum\limits_{k-1}^mu_k,
\]
then we define
\[
Au=\sum\limits_{k-1}^mA^{(k)}u_k
\]
This operator $A$ will be a generating operator.
$\triangle$

Such operator $A$ is called an elliptic operator if the operator-function $A(x)$ consists of Fredholm operators $\forall x\in M$.
In a certain sense we can obtain the inverse result.

\section{The index theorem}

Here we will give an index theorem for our operators. It seems it does not give real instrument for calculating index, but it shows us what kinds of operators we need to study for obtaining good index formulas.

{\bf Theorem 2.} {\it
The index of the operator $A$ on the manifold $M$ is a sum of corresponding indices
\begin{equation}\label{3}
Ind~ A=\sum\limits_{k=0}^{m}Ind~A^{(k)},
\end{equation}
}

{\bf Proof.}
Indeed, all operators $A^{(k)}, k=0,1,\cdots,m$ act in different spaces. Therefore, the generating operator $A$ has the following kernel and co-kernel
\[
Ker~A=\sum\limits_{k=0}^{m}Ker~A^{(k)}
\]
\[
Coker~A=\sum\limits_{k=0}^{m}Coker~A^{(k)}
\]
According to definition for an index we obtain the formula \eqref{3}.
$\triangle$

\section{Example: pseudo-differential constructions}

\subsection{Local situations}

We consider a certain integro-differential operator $A$ on $m$-dimensional compact manifold $M$ with a boundary. This operators is defined by the function $A(x,\xi), (x,\xi)\in{\bf R}^{2m}$. We will suppose that the symbol has the order $\alpha\in{\bf R}$, i.e.
\begin{equation}\label{4}
c_1(1+|\xi|)^{\alpha}\leq|A(x,\xi)|\leq c_2(1+|\xi|)^{\alpha},
\end{equation}
for all admissible $x,\xi$ with universal positive constants $c_1,c_2$.

We consider such a compact manifold $M$ with a boundary that there are some smooth compact sub-manifolds  $M_k$ of dimension $0\leq k\leq m-1$ on the boundary  $\partial M$ of manifold  $M$ which are singularities of a boundary. These singularities are described by a local representative of operator $A$ in a point  $x_0\in M$ on the map $U\ni x_0$ in the following way
\begin{equation}\label{5}
(A_{x_0}u)(x)=\int\limits_{D_{x_0}}\int\limits_{{\bf R}^m}e^{i\xi\cdot(x-y)}A(\varphi(x_0),\xi)u(y)d\xi dy,~~~x\in D_{x_0},
\end{equation}
where $\varphi :U\to D_{x_0}$ is a diffeomorphism, and the canonical domain $D_{x_0}$ has a distinct form depending on a placement of the point $x_0$ on manifold  $M$. We consider following canonical domains  $D_{x_0}$: ${\bf R}^m, {\bf R}^m_+=\{x\in{\bf R}^m: x=(x',x_m), x_m>0\}, W^k={\bf R}^k\times C^{m-k}$, where $C^{m-k}$ is a convex cone in  ${\bf R}^{m-k}$. For instance, if we consider a cube $Q$ in 3-dimensional space then we have 4 canonical domains: ${\bf R}^3$ for inner points, ${\bf R}^3_+=\{x\in{\bf R}^3: x=(x_1,x_2,x_3), x_3>0\}$ for six 2-faces, ${\bf R}\times C^2=\{x\in{\bf R}^3: x=(x_1,x_2,x_3), x_>0,x_3>0\}$ for twelve 1-dimensional edges, and $C^3=\{x\in{\bf R}^3: x=(x_1,x_2,x_3), x_1>0,x_2>0,x_3>0\}$ for eight vertices.

Such an operator $A$ will be considered in Sobolev--Slobodetskii spaces $H^s(M)$, and local variants of such spaces will be spaces  $H^s(D_{x_0})$.

{\bf Definition 3.} {\it
The symbol of an operator $A$ is called the operator-function  $A(x): M\rightarrow\{A_x\}_{x\in{M}}$ which is defined by local representatives of the operator $A$.
}

Under some additional assumptions on smoothness properties of the function $A(x,\xi)$ one has the following

{\bf Theorem 3.} {\it
The operator  $A$  has a Fredholm property iff its symbol is composed by Fredholm operators.
}

Simplest variant of this theorem was proved in \cite{S,V1}. For general local operators in Lebesgue spaces Theorem 3 was proved in \cite{S1}.

{\bf Definition 4.} {\it
An operator $A$ is called an elliptic operator if its symbol is composed by invertible operators.
}

{\bf Corollary 1.} {Elliptic operator is a Fredholm operator.
}

{\bf Remark 3.} {\it
~If an ellipticity property does not hold on sub-manifolds  $M_k$ one needs to modify local representatives of the operator $A$ adding special
boundary or co-boundary operators.
}

Using a special partition of a unity on the manifold  $M$, elliptic symbol $A(x)$ for each $x\in\overline{M_k}$ which is given ny the formula \eqref{5} and above constructions from Theorem 1 we
 obtain $m+1$ operators
$A^{(k)}$ according to a number of singular sub-manifolds including whole boundary  $\partial M$ and the manifold $M$.

{\bf Theorem 4.} {\it
Index of the Fredholm pseudo-differential operator $A$  is given by the formula
\[
Ind~A=\sum\limits_{k=0}^{m}Ind~A^{(k)}.
\]
}

{\bf Proof.}
Really, this is a simply corollary from Theorem 2. Indeed, we need to show exactly what spaces we choose as $H_k^{(j)}, j=1,2; k=0,1,\cdots,m$. We enumerate:
\[
A^{(m)}: H^s({\bf R}^m)\rightarrow  H^{s-\alpha}({\bf R}^m);
\]
\[
A^{(m-1)}: H^s({\bf R}^m_+)\rightarrow  H^{s-\alpha}({\bf R}^m_+)'
\]
\[
A^{(k)}: H^s(W^k)\rightarrow  H^{s-\alpha}(W^k), k=0,1,\cdots m-2,
\]
so that $H_m^{(1)}= H^s({\bf R}^m), H_m^{(2)}=H^{s-\alpha}({\bf R}^m), H_{m-1}^{(1)}= H^s({\bf R}^m_+), H_{m-1}^{(2)}=H^{s-\alpha}({\bf R}^m_+), H_k^{(1)}= H^s(W^k), H_k^{(2)}=H^{s-\alpha}(W^k)$. Then we compose the direct sum of such spaces and the operator $A'$ acting in these direct sums
\[
A': H^s({\bf R}^m)\oplus H^s({\bf R}^m_+)\oplus\sum\limits_{k=0}^{m-2} H^s(W^k)\longrightarrow 
\]
\[
H^{s-\alpha}({\bf R}^m)\oplus H^{s-\alpha}({\bf R}^m_+)\oplus\sum\limits_{k=0}^{m-2} H^{s-\alpha}(W^k)
\]
Let us note that the operator $A'$ doesn't coincide with the operator $A$, but these operators have the same local representatives, i.e. the same symbols. We call the operator $A'$ {\it virtual representative} of the operator $A$. Since homotopies of symbols one-to-one correspond to homotopies of operators we complete the index theorem.
$\triangle$

Of course Theorem 4 does not give effective index formulas, but it shows what kinds of operators we need to consider from index theory viewpoint.

{\bf Remark 4.} {\it
~If we consider an elliptic pseudo-differential operator in $H^s({\bf R}^m_+)$ \cite{E} with the smooth symbol $A(x,\xi)$ we have two decomposition operators: $A^{(m)}$ related to closure of inner points of ${\bf R}^m_+$ and $A^{(m-1)}$ related to boundary points ${\bf R}^{m-1}$. Operator symbols are distinct nature for inner and boundary points. For the first case such a symbol is represented by integral over the whole ${\bf R}^m$, but for the second case this integral is taken for a half-space. The index of $A^{(m)}$ will be zero according to classical Atiyah--Singer theorem, but the index of $A^{(m-1)}$ depends on so-called index of factorization for the symbol $A(x,\xi)$ at boundary point $x\in{\bf R}^{m-1}$.
}

\subsection{The wave factorization: harmonic analysis and complex variables}

To obtain invertibility conditions for local operators we need some additional characteristics for the classical symbol of elliptic pseudo-differential operators. The studying invertibility of a local operator in $W^k$ or in other words the unique solvability of the equation
\[
(A_{x_0}u)(x)=v(x),~~~x\in W^k,
 \]
in Sobolev--Slobodetskii space $H^s(W^k)$   is equivalent to a unique solvability for so called paired equation
\begin{equation}\label{6}
(A_{x_0}P_+U)(x)+(IP_-U)(x)=V(x),~~~x\in{\bf R}^m,
\end{equation}
in the space $H^s({\bf R}^m)$, where $P_+, P_-$ are projectors on $W^k, {\bf R}^m\setminus W^k$, it can be easily proved. And now if we apply the Fourier transform then we will come to complex spaces \cite{V1}

We denote by $\stackrel{*} {C^{m-k}}$ the conjugate cone for the $C^{m-k}$:
\[
\stackrel{*} {C^{m-k}}=\{x\in{{\bf R}}^m: x\cdot y>0, \forall y\in C^{m-k}\},
\]
$T(\pm\stackrel{*} {C^{m-k}})$ denotes a radial tube domains over the cone $\pm\stackrel{*} {C^{m-k}}$ \cite{Vl}, i.e. a domain of multidimensional complex space ${{\bf C}}^m$ of the type ${{\bf R}}^m\pm\stackrel{*} {C^{m-k}}$.

Let the classical symbol $a(\xi), \xi\in{\bf R}^m,$ in local coordinates satisfies the condition
\[
c_1(1+|\xi|)^{\alpha}\leq|a(\xi)|\leq c_1(1+|\xi|)^{\alpha}.
\]
Let us denote $\xi=(\xi'',\xi'), \xi''=(\xi_1,\cdots,\xi_k), \xi'=(\xi_{k+1},\cdots,\xi_m)$.

{\bf Definition 5.} {\it
$k$-wave factorization of the symbol $a(\xi)$ with respect to the cone  $C^{m-k}$ is called its representation in the form
$$
a(\xi)=a_{\neq}(\xi)a_=(\xi),
$$
where the factors  $a_{\neq}(\xi),a_=(\xi)$ must have the following properties:

1) $a_{\neq}(\xi), a_=(\xi)$ are defined for all  $\xi\in{\bf R}^{m}$ excluding may be the points  ${\bf R}^{k}\times\partial\left(\stackrel{*} {C^{m-k}}\cup( -\stackrel{*} {C^{m-k}})\right)$;

2) $a_{\neq}(\xi), a_=(\xi)$ admit analytical continuation into radial tube domains  $T(\stackrel{*} {C^{m-k}}),T(-\stackrel{*} {C^{m-k}})$ for almost all  $\xi''\in{{\bf R}}^k$ respectively with estimates
$$
|a_{\neq}^{\pm 1}(\xi'',\xi'+i\tau)|\leq c_1(1+|\xi|+|\tau|)^{\pm\ae_k},
$$
$$
|a_{=}^{\pm 1}(\xi'',\xi'-i\tau)|\leq c_2(1+|\xi|+|\tau|)^{\pm(\alpha-\ae_k)},~\forall\tau\in\stackrel{*} {C^{m-k}}.
$$

The number $\ae_k\in{{\bf R}}$ is called an index of  $k$-wave factorization.
}

\subsection{Fredholm properties}

For simplicity we consider here the case when $M$ is a bounded domain in ${\bf R}^m$ and its classical symbol looks as $A(x,\xi)$.
Here we assume additionally that symbol of the operator $A$ is continuous on $M_k, k=0,1,...,m,$  family of operators (of course with respect to the norm $|||\cdot|||$). This property holds for example if the function  $A(x,\xi), (x,\xi)\in M\times{\bf R}^m$ is continuous differentiable up to boundary. Then according to enveloping theorem \cite{S} using operator symbol one can construct  $n$ operators $A_k$. If these operators have Fredholm property then the general operator will have a Fredholm property with the index according to Theorem 4.

Let $\ae_{n-1}(x)$ be the index of factorization \cite{E} of the function $A(x,\xi)$ in the point $x\in\partial M\setminus\cup_{k=0}^{m-2}M_k$,
$\ae_k(x)$ be indices of $k$-wave factorization with respect to the cone $C^{m-k}_x$ at points  $x\in M_k, k=0,1,\cdots,m-2$
and we assume that the functions  $\ae_k(x), k=0,1,\cdots,m-1,$ are continuously continued in  $\overline{M_k}$.

{\bf Remark 5.} {\it
Similarly \cite{E} using uniqueness result for the wave factorization \cite{V1} one can verify that the functions $\ae_k(x), k=0,1,\cdots,m-1,$ don't depend on local coordinates
}

{\bf Theorem 5.} {\it
If the classical elliptic symbol  $A(x,\xi)$ admits $k$-wave factorization with respect to the cones $C^{m-k}$ with indices  $\ae_k(x), k=0,1,\cdots,m-2$ satisfying the condition
\begin{equation}\label{7}
|\ae_k(x)-s|<1/2,~~~\forall x\in M_k,~~~k=0,1,\cdots,m-1,
\end{equation}
then the operator $A: H^s(M)\rightarrow H^{s-\alpha}(M)$ has a Fredholm property.
}

{\bf Proof.}
To prove the theorem we need to verify invertibility properties for all local representatives for our pseudo-differential operator $A$.

\underline{\it A whole space.} This case was historically the first in the theory of pseudo-differential equations. If $x_0\in\stackrel{\circ}{M}$
is an inner point then the local representative in the formula \eqref{5} has the following form (in local coordinates $\varphi$)
\[
(A_{x_0}u)(x)=\int\limits_{{\bf R}^m}\int\limits_{{\bf R}^m}e^{i\xi\cdot(x-y)}A(\varphi(x_0),\xi)u(y)d\xi dy,~~~x\in{\bf R}^m,
\]
and this is classical pseudo-differential operator.\cite{E,RS,V1}. Ellipticity condition for the classical symbol
$$
A(x,\xi)\neq 0
$$
 for all admissible $x,\xi$ is necessary and sufficient condition for invertibility of every such operator.

 Unfortunately, if we have a piece of the space ${\bf R}^m$ we need to study invertibility properties for the operator in left-hand side of the equation \eqref{6}.

\underline{\it A half-space.} If $x_0\in\partial M$ is smoothness point of $\partial M$ then a local representative for the operator $A$ has the following form
\[
(A_{x_0}u)(x)=\int\limits_{{\bf R}^m_+}\int\limits_{{\bf R}^m}e^{i\xi\cdot(x-y)}A(\varphi(x_0),\xi)u(y)d\xi dy,~~~x\in{\bf R}^m_+,
\]
To study solvability for corresponding paired equation \eqref{6} there were used a factorization theory and one dimensional singular integral operators \cite{E,G,M}. Full solvability theory for such equations was constructed in M.I. Vishik--G.I. Eskin papers (see \cite{E}). Principal role takes the index of factorization, in our notation $\ae_{m-1}$, if the condition \eqref{7} holds then the operator$H^s({\bf R}^m_+)\rightarrow H^{s-\alpha}({\bf R}^m_+)$ is invertible.

\underline{\it A cone.} Here we have more complicated local representative
\[
(A_{x_0}u)(x)=\int\limits_{W^k}\int\limits_{{\bf R}^m}e^{i\xi\cdot(x-y)}A(\varphi(x_0),\xi)u(y)d\xi dy,~~~x\in W^k,
\]
but the factorization idea works here also in multidimensional context; if $k$-wave factorization exists then the condition \eqref{7} is sufficient for invertibility of such operator \cite{V1}.
$\triangle$

{\bf Remark 6.} {\it
If the ellipticity does not hold on sub-manifolds  $M_k$ then we can modify the operator $A$ using boundary or co-boundary operators \cite{V1}. Particularly we need such constructions if one of conditions $\eqref{7}$ does not hold.
}

Some considerations related to this paper are given in \cite{V2,V3,V4,V5}, particularly these are related to more complicated singularities and more general spaces.

\section{Conclusion}

We have described a new approach to constructing the theory of pseudo-differential equations and related boundary value problems. This approach is based on general principles for special local operators.

In our opinion such considerations will be useful for discrete situations in which pseudo-differential operators are defined in functional spaces of discrete variable. Some first considerations in this direction were done, for example, in \cite{VV1}. Moreover, discrete situation is more accessible, since it  permits to apply computer calculations. We hope to develop these studies in this direction including a comparison between discrete and continuous cases.



\end{document}